\documentclass[11pt]{amsart} 

\usepackage[all]{xy}
\usepackage{amsthm}
\usepackage{mathrsfs}
\usepackage{amssymb}
\usepackage{amsmath} 
\usepackage{array, url}
\usepackage{stmaryrd}
\usepackage{mathtools}
\usepackage{tensor}
\usepackage{stmaryrd}
\usepackage{pifont}
\usepackage{float}
\usepackage[labelformat=empty]{caption}
\usepackage[margin=1in]{geometry}

\let\oldmarginpar\marginpar
\renewcommand\marginpar[1]{\-\oldmarginpar[\raggedleft\footnotesize #1]%
{\raggedright\footnotesize #1}}

\newcounter{firstnumber}[section]
\newcounter{secondnumber}[firstnumber]
\newcounter{thirdnumber}[secondnumber]
\newcounter{fourthnumber}[thirdnumber]
\newcounter{fifthnumber}[fourthnumber]
\newcounter{currentdepth}

\renewcommand{\thefirstnumber}{\arabic{section}.\arabic{firstnumber}}

\newcommand{\segment}[2]{%
    \setcounter{currentdepth}{1}%
    \def\thesubsection{\thefirstnumber}%
    \refstepcounter{firstnumber}\label{#1}
    \addtocounter{subsection}{-1}%
    \subsection{#2}}

\swapnumbers

\newcommand{\Spec}{\operatorname{Spec}}

\newcommand{\Pic}{\operatorname{Pic}}
\newcommand{\Div}{\operatorname{Div}}

\newcommand{\Hom}{\operatorname{Hom}}

\newcommand{\m}[1]{\mathrm{#1}}

\newcommand{\bb}[1]{\mathbb{#1}}

\newcommand{\ka}{\kappa}
\newcommand{\si}{\sigma}
\newcommand{\ze}{\zeta}
\newcommand{\ga}{\gamma}
\newcommand{\al}{\alpha}
\newcommand{\be}{\beta}

\newcommand{\Gm}{{\mathbb{G}_m}}
\newcommand{\Qp}{{\QQ_p}}
\newcommand{\Zp}{{\ZZ_p}}
\newcommand{\Fp}{{\FF_p}}
\newcommand{\Qbar}{{\overline{\QQ}}}

\newcommand{\ZZ}{\bb{Z}}

\newcommand{\QQ}{\bb{Q}}
\newcommand{\PP}{\bb{P}}

\newcommand{\FF}{\bb{F}}
\renewcommand{\AA}{\bb{A}}

\newcommand{\Oo}{\mathcal{O}}
\newcommand{\Kk}{\mathcal{K}}

\newcommand{\SSS}{\mathbf{S}}

\newcommand{\inv}{^{-1}}

\newcommand{\zpn}{\ze_{p^n}}
\newcommand{\Qzpn}{\QQ(\zpn)}

\newcommand{\HTn}{G_{T,n}}

\newcommand{\HTinfty}{G_{T,\infty}}

\newcommand{\Zpn}{{\ZZ/p^n}}

\newcommand{\Zpstar}{\ZZ_p^*}
\newcommand{\Zpnstar}{(\ZZ/p^n)^*}
\newcommand{\Fpstar}{\FF_p^*}

\newcommand{\mupn}{\mu_{p^n}}

\newcommand{\Ttilde}{{\widetilde{T}}}

\newcommand{\et}{{\textrm {\'et}}}

\begin{document}

\SelectTips{cm}{11}

\title[The Heisenberg coboundary equation]{The Heisenberg coboundary equation: appendix to \textit{Explicit Chabauty-Kim theory}}
\author{Ishai Dan-Cohen and Stefan Wewers}
\date{\today}
\maketitle

\begin{abstract}

Let $p$ be a regular prime number, let $G_{\{p\}}$ denote the Galois group of the maximal unramified away from $p$ extension of $\QQ$, and let $H_\et$ denote the Heisenberg group over $\Qp$ with $G_{\{p\}}$-action given by $H_\et = \QQ_p(1)^2 \oplus \QQ_p(2)$. Although Soul\'e vanishing guarantees that the map $ H^1(G_{\{p\}}, H_\et) \to  H^1(G_{\{p\}}, \Qp(1)^2)  $ is bijective, the problem of constructing an explicit lifting of an arbitrary cocycle in $H^1(G_{\{p\}}, \Qp(1)^2)  $ proves to be a challenge. We explain how we believe this problem should be analyzed, following an unpublished note by Romyar Sharifi, hereby making the original appendix to \textit{Explicit Chabauty-Kim theory} available online in an arXiv-only note.

\end{abstract}

\section{The context}

This brief note began its life as an appendix to \textit{Explicit Chabauty-Kim theory for the thrice punctured line in depth two} \cite{CKTwo}, which received an appendectomy prior to publication. Let $X = \PP^1\setminus\{0,1,\infty\}$, let
\[
S = \{q_1, \dots, q_s\}
\]
denote a finite set of primes, let
\[
\SSS = \Spec \ZZ \setminus S
\,,
\]
let $p$ denote a prime $\notin S$, and let $T = S \cup \{p\}$. Kim's approach to the study of the set $X(\SSS)$ of $S$-integral points of $X$ involves a certain tower of morphisms of affine finite-type $\Qp$-varieties. As we explain in \textit{Explicit Chabauty-Kim theory}, its first two steps look like so. 
\[
\xymatrix{
H^1_f(G_T, H_\et) \ar[r]^-{h_2} \ar[d]_\cong^{\pi_*} & \AA^3_\Qp \ar[d] \\
H^1_f(G_T, \QQ_p(1)^2)\ar[r]_-{h_1} & \AA^2_\Qp
}
\]
Here $G_T$ denotes the Galois group of the maximal unramified outside of $T$ extension of $\QQ$, $H_\et$ denotes the Heisenberg group object with $G_T$-action given simply by $H_\et = \Qp(1)^2 \oplus \Qp(2)$, and the $H^1_f$'s are certain subschemes of nonabelian cohomology varieties. The map induced by abelianization
\[
\pi:H_\et \to \Qp(1)^2
\]
on the level of $H^1_f$'s fits in a commuting square like so.
\[
\xymatrix{
H^1_f(G_T, H_\et) \ar[d]_\cong^{\pi_*} \ar@{^(->}[r] &
H^1(G_T, H_\et) \ar[d]_\cong^{\pi_*} \\
H^1_f(G_T, \Qp(1)^2) \ar@{^(->}[r] \ar@{=}[d] &
H^1(G_T, \Qp(1)^2) \ar@{=}[d] \\
\QQ_p^S  \ar@{^(->}[r] &
\QQ_p^T
}
\]
Our main result in \textit{Explicit Chabauty-Kim theory} is a complete computation of the map
\[
h_2 \circ \pi_*\inv : \QQ_p^S \to \QQ_p^3
\,.
\]
Since our methods there where somewhat indirect, we document here our initial attempt to compute $\pi_*\inv$ directly.

\section{The Problem}

\segment{110814c}{}
By Corollary 6.2.3 of \cite{CKTwo}, an $A$-point of $\m H^1_{ \m f} (G_T, \Qp(1))$ may be written $\ka_x$, where 
$$
x = q_1^{x_1}\cdots q_s^{x_s}
$$
is a formal product of powers, with $x_1, \dots, x_r \in A$. If 
$$
y = q_1^{y_1}\cdots q_s^{y_s}
$$
denotes another point, the cup product $\ka_y \cup \ka_x$ is an element of $\m Z^1(G_T, A(2))$. For simplicity, restrict attention to the case $A= \Qp$, and consider the cochain complex
$$
0 \to \m C^0(\Qp(2)) \to \m C^1(\Qp(2)) \to \m C^2(\Qp(2)) \to \m C^3(\Qp(2)) \to \cdots \;.
$$
for the cohomology of $G_T$ with coefficients in $\Qp(2)$. By the vanishing results
$$
\m H^1(G_T, \Qp(2)) = \m H^2(G_T, \Qp(2)) = 0 \,,
$$
the equation
\[
\ka_y \cup \ka_x = d\al
\]
in $\m C^2$ admits a solution $\al \in \m C^1$, unique up to translation by the coboundary
$
d\beta: \si \mapsto \beta - \si(\beta)
$
of an element $\beta \in \Qp(2)$. Recalling the definition of the cup product and the second coboundary, we have
\[
\ka_x(\si) \otimes \si\ka_y(\tau) = \al(\si\tau) - \si\al(\tau) - \al(\si)
\,.
\]
Here, $\si$ and $\tau$ vary over $G_T$. We call this \textit{the Heisenberg coboundary equation}.

\segment{120130-1}{}
By \S5 of \cite{CKTwo}, we have an exact sequence
\[
\xymatrix{
1 		\ar[r]					&
\Qp(2) 	\ar[r] 					&
H_\et 	\ar[r]					&
\Qp(1)^2	\ar[r] \ar@/_/[l]_\Sigma		&
1								}
\]
of nonabelian $G_T$-modules, which is split if we identify $H_\et$ with its Lie algebra and forget the bracket. Fix a point $(x,y) \in (\ZZ[S\inv]^*\otimes \Qp)^2$ in the source of the unipotent $p$-adic Hodge morphism in depth one. We then have associated Kummer cocycles $\ka_x, \ka_y: G_T \rightrightarrows \Qp(1)$, and by composing with $\Sigma$, we obtain a candidate $\Sigma(\ka_x, \ka_y) \in \m C^1(G_T, H_\et)$ for a lifting of $(\ka_x, \ka_y)$ to depth two. By segment 2.3.2 of \cite{CKTwo}, its failure to be a cocycle is measured by a solution $\al$ of the Heisenberg coboundary equation. So if we set $\ka_{x,y} := \al\inv \cdot \Sigma (\ka_x, \ka_y)$, we obtain a representative for the element of $\m H^1(G_T, H_\et)$ which maps to $(\ka_x, \ka_y)$ in $\m H^1(G_T, \Qp(1)^2)$.

\segment{120130-2}{}
The problem then, is to make the solution $\al$ of the Heisenberg coboundary equation in some way explicit.

\section{Steps towards its solution}

\segment{120130-3}{}
Soul\'e's proof of the vanishing of $\m H^2(\Qp(2))$ is not well adapted to this application. A simpler proof is given by Romyar Sharifi in an unpublished note \cite{Sharifi}, for the case $p=2$. As Sharifi points out, the essential property of the even prime which makes his proof possible is its regularity. Sharifi's proof goes roughly as follows. Let $K_T$ denote the maximal unramified outside $T$ extension of $\QQ$. Let $\HTn$ denote the Galois group of $K_T$ over $\Qzpn$. Similarly, we let $\HTinfty$ denote the Galois group of $K_T$ over $\QQ(\ze_{p^\infty})$. Noting that the Galois group of $\Qzpn/\QQ$ is $(\ZZ/p^n)^*$, we have the following tower of fields and Galois groups for each $n$. 
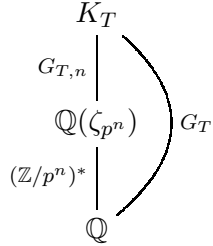
\[\xymatrix{
K_T \ar@{-}[d]_\HTn	\ar @{-} @/^6ex/ [dd]^{G_T}		\\
\Qzpn \ar@{-}[d]	_{(\ZZ/p^n)^*}		\\
\QQ				}
\]
By direct computation applied to the low degree terms of the Hochschild-Serre spectral sequence, we obtain an isomorphism
\[
\m H^1(G_T, \Qp(2)) = \m H^1(\HTinfty, \Qp(2))^{\Zpstar} \tag{$\star$} \;.
\]
On the other hand, for $p$ regular, we have
\[
\Qp \otimes_\Zp (\ZZ[T\inv, \ze_{p^\infty}]^*_{/p})(1) = \m H^1(\HTinfty, \Qp(2)) \tag{$\star\star$}  \;.
\]
The subscript $/p$ indicates $p$-adic completion. The argument here may be summarized as follows: there's always an injection from the left to the right coming from the Kummer exact sequence; the cokernel lives inside the Picard group (suitably interpreted), whose (pro-)order is (pro-)coprime to $p$. A study of the action of $\ZZ_p^*$ on $\Qp \otimes_\Zp (\ZZ[T\inv, \zpn]^*_{/p})(1)$ now leads to the conclusion that
$$
\m H^1(G_T, \Qp(2)) = \Qp \,.
$$
Finally, Poitou-Tate duality is used as a vehicle to get to $\m H^2$.

Actually, throughout most of the proof, Sharifi works with $n$ finite. For $n$ finite, statements analogous to ($\star$), ($\star\star$) fail. Their failure however, is measured by groups whose order turns out to be finite and bounded in $n$.

\segment{110522c}{}
Sharifi's use of Poitou-Tate duality presents for us an obstacle. On the other hand, since \textit{many} regular primes are known to exist (see, for instance \S5.3 of Washington \cite{Washington}), the stipulation that $p$ be regular is relatively harmless. So a possible approach may be to attack the vanishing of $\m H^2$ (or at least of the relevant elements of $\m H^2$) directly, by methods inspired by Sharifi's computation of $\m H^1$. To do so, we would replace \ref{120130-3}($\star$) by an analysis of the map
\[
\m H^2(G_T, \Qp(2)) \to \m H^2(G_{T \infty}, \Qp(2))^{\Zpstar}
\tag{\ding{65}}
\,,
\]
and we would replace \ref{120130-3}($\star \star$) by the map
\[
\Qp\otimes_{\Zp} K_2^\m{M} (\ZZ[\ze_{p^\infty}, T\inv ]) \to \m H^2(G_{T\infty}, \Qp(2))
\tag{\ding{65}\ding{65}}
\,,
\]
while keeping track of the $\Zpstar$ action.

\section{An ensuing family of spectral sequences in Galois cohomology}

\segment{110502f}{}
If \ref{110522c}(\ding{65}) fails to be bijective, the failure is best measured by certain terms in an associated family of spectral sequences. Let
$$
1 \to N \to G \to Q \to 1
$$
be a short exact sequence of (topological) groups, and $A$ a (continuous) $\ZZ[G]$-module whose addition law we denote by $\star$. Then the (continuous) cohomology groups $\m H^q(N, A)$ have a natural structure of (continuous) $\ZZ[Q]$-module, and there's a spectral sequence
$$
E_2^{p,q}= \m H^p(Q, \m H^q(N, A)) \Rightarrow \m H^{p+q}(G, A) \;.
$$
Elements of $\m H^1(N, A)$ may be represented by (continuous) maps $\phi:N\to A$ which satisfy 
$$
\phi(\si\tau) = \phi(\si)\star\si\phi(\tau) \;.
$$
If $\phi$ is such a map and $\al$ is an arbitrary element of $Q$, then the action of $Q$ on $\m H^1(N, A)$ is given in terms of cocycles by lifting $\alpha$ arbitrarily to an element $\ga$ of $G$ and declaring that for any $\eta \in N$,
$$
\phi^\al(\eta) = \alpha\inv\phi(\ga\eta\ga\inv) \;.
$$
See \S5, 6 of Chapter VII of \cite{corloc}.

\segment{120131-1}{}
For each $n>m$, we may apply this to the short exact sequence
\[
1 \to G_{T,n} \to G_T \to \Zpnstar \to 1
\,,
\]
with coefficients in $(\ZZ/p^m)(2)$. If we set
\[
\tensor*[^m_n]{E}{^{p,q}_2} := \m H^p(\Zpnstar; \m H^q(G_{T, n}; \ZZ/p^m(2)))
\]
and
\[
\tensor*[^m_n]{H}{^r} := \m H^r(G_T, \ZZ/p^m(2))
\,,
\]
then there's a spectral sequence
\[
\tensor*[^m_n]{E}{^{p,q}_2} \Rightarrow \tensor*[^m_n]{H}{^r}
\,.
\]
Relevant terms and arrows of this spectral sequence are pictured below.
\begin{small}
\[\xymatrix @ R=25pt @ C=0pt {
\tensor*[^m_n]{E}{^{0,2}_2} \ar[drr]		&
						&
						&
						&
\tensor*[^m_n]{E}{^{0,2}_3} \ar[ddrrr]	&
						&
						&
						\\
						&
\tensor*[^m_n]{E}{^{1,1}_2}	\ar[drr]		&
\tensor*[^m_n]{E}{^{2,1}_2}				&
						&
						&
						&
						&
						\\
						&
						&
\tensor*[^m_n]{E}{^{2,0}_2}				&
\tensor*[^m_n]{E}{^{3,0}_2}				&
						&
						&
						&
\tensor*[^m_n]{E}{^{3,0}_3}				}
\]
\end{small}

\segment{110429b}{}
We now discuss the terms $\tensor*[^m_n]{E}{^{2,0}_2}$. We set $n=m$ for simplicity.

\subsection*{Proposition}
Each $^n_nE_2^{2,0}=H^2(\Zpnstar, H^0(\HTn, \Zpn(2)))$ is a finite group of bounded order.

\medskip
\noindent
The proof is in segments \ref{110429c}--\ref{110502c}.

\segment{110429c}{}
$\HTn$ acts trivially on $\Zpn(2)$, so
$$
\m H^0(\HTn, \Zpn(2)) = \Zpn(2) \,.
$$
We have a short exact sequence of groups
$$
0 \to 1+(p) \to \Zpnstar \to \Fpstar \to 0
$$
from which we obtain a spectral sequence
$$
F_2^{p,q}=\m H^p(\Fpstar, \m H^q(1+(p), \Zpn(2))) \Rightarrow \m H^{p+q}(\Zpnstar, \Zpn(2)).
$$
It suffices to show that the terms $F_2^{2,0}$, $F_2^{1,1}$, $F_2^{0,2}$ are finite groups of bounded order. But since $\Fpstar$ itself is finite cyclic of bounded order, it suffices to show that the three cohomologies $\m H^0, \m H^1, \m H^2(1+(p), (\Zpn(2)))$ are finite groups of bounded order.

\segment{110502a}{}
Let $C$ be a finite cyclic group with generator $\si$, consider the elements $1-\si$, $N:= \sum_{\tau\in C} \tau$ of the group algebra $\ZZ[C]$, and let $A$ be a $\ZZ[C]$-module. Then the sequence
$$
0 \to A \xrightarrow{\si-1} A \xrightarrow{N} A \xrightarrow{\si-1} A \xrightarrow{N} \cdots \;,
$$
in which the first $A$ is in degree zero, forms a complex $A^\bullet$ and
$$
\m H^i (C,A)= \m H^iA^\bullet \;.
$$

\segment{110502c}{}
Returning to the situation and the notation of the proposition, we note that $1+(p)$ is generated by the element $e^p = 1+p+\frac{p^2}{2!}+\cdots$ and that $e^p$ acts on $\Zpn(2)$ by multiplication by $e^{2p}$. Thus, to complete the proof of the proposition, we need only note that (under our assumption that $p\neq2$)
$$
v_p(e^{2p} -1) = 1 \;,
$$
so that the endomorphism of $\Zpn$ given by multiplication by $e^{2p}-1$ has kernel $(p^{n-1})$ and cokernel $\Fp$, both of which have order $p$, hence in particular bounded, as hoped.

\segment{120131-3}{}
For the remainder of the section we focus our attention on the terms
\[
\tensor*[^m_n]{E}{^{1,1}_2} = \m H^1(\Zpnstar; \m H^1(G_{T,n};\ZZ/p^m(2)))
\,.
\]
We again set $n=m$ for simplicity.

\segment{110502g}{}
We denote the group $\mupn(\QQ(\zpn))$ of $(p^n)^\m{th}$ roots of $1$ in $\QQ(\zpn)$ by $\mupn$ for short. $\mupn$ is isomorphic to $\Zpn(2)$ as a $\Zpn$-module. 
Moreover, if we let an arbitrary element $\al$ of $\Zpnstar$ act on an arbitrary element $\ze$ of $\mupn$ by $\ze^{\al^2}$, then any such isomorphism becomes equivariant with the action of $\Zpnstar$. Recalling our formula for the action of a quotient group on the first cohomology of the kernel in terms of cocycles (\ref{110502f}), and noting that $\HTn$ acts trivially on $\Zpn(2)$, we obtain an isomorphism
$$
\m H^1(\HTn, \Zpn(2)) \cong \Hom(\HTn, \mupn)
$$
which is not canonical, but is nevertheless equivariant for the action of $\Zpnstar$ on $\Hom(\HTn, \mupn)$ given in terms of a continuous map
$$
\phi: \HTn \to \mupn \;,
$$
an $\eta \in \HTn$, an $\al \in \Zpnstar$, and a lifting $\ga$ of $\al$ to $G_T$, by the formula
$$
\phi^\al(\eta) = (\phi(\ga \eta \ga \inv))^{\al^{-2}} \,.
$$

\segment{110502h}{}
Let 
$$
E= \{ a \in \QQ(\zpn)^* \; | \; v(a) \equiv 0 \mod p^n \hspace{3mm} \forall v \nmid T  \} \;.
$$
Given an element $\al \in \Zpnstar$, and an element $a \in E/\QQ(\zpn)^{*p^n}$, we let $\al$ act on $a$ by
$$
\al\inv(a)^{\al\inv} \;.
$$
Here, the $\al\inv$ on the left acts on $a$ through the Galois action of $\Zpnstar$ on $\QQ(\zpn)$, while the $\al\inv$ in the exponent (which may equivalently be put inside the parentheses) denotes multiplication of the base by itself ``$\al\inv$ many times", an operation which is only well defined modulo $\QQ(\zpn)^{*p^n}$. 

\segment{110502i}{}
An element $a \in E$ gives rise to a Kummer cocycle $\ka_a$, which is unramified outside $T$. This means that $\ka_a$ defines a map $\HTn \to \mupn$ given in terms of an element $\eta \in \HTn$ and a $(p^n)^\m{th}$ root $a^{1/p^n}$ of $a$, by the formula
$$
\ka_a(\eta) = \frac {\eta (a^{1/p^n})} {a^{1/p^n}} \;.
$$

\segment{110502e}{Proposition}
In the notation and the situation of paragraphs \ref{110502h} and \ref{110502i}, the assignment
$$
a \mapsto \ka_a
$$
defines a $\Zpnstar$-equivariant isomorphism
$$
E/\QQ(\zpn)^{*p^n} \xrightarrow{\cong} \Hom(\HTn, \mupn) \;.
$$

\medskip
\noindent
The proof is in segments \ref{110503c}--\ref{110503b}.

\segment{110503c}{}
Fix an algebraic closure $\Qbar$ of $\QQ(\zpn)$, let $H$ denote the Galois group of $\Qbar$ over $\QQ(\zpn)$ and let $N$ denote the Galois group of $K_T$ (\ref{120130-3}) over $\QQ(\zpn)$:
\[\xymatrix{
\Qbar \ar@{-}[d]_N	\ar @{-} @/^4ex/ [dd]^H	\\
K_T \ar@{-}[d]	_\HTn		\\
 \QQ(\zpn) 			}
\]
Evaluating the Kummer exact sequence
$$
1 \to \mupn \to \Gm \to \Gm \to 1
$$
on $\Qbar$, applying invariants with respect to the action of $H$, recalling Hilbert's theorem 90, and noting that $H$ acts trivially on $\mupn$, we obtain an isomorphism
$$
\ka:\QQ(\zpn)^*/\QQ(\zpn)^{*p^n} \xrightarrow{\cong} \Hom(H, \mupn) \;.
$$
Then
$$
\ka\inv(\Hom(\HTn, \mupn))=E/\QQ(\zpn)^{*p^n} \;.
$$
Indeed, given $a \in E$, $\ka_a$ factors through $\HTn$ if and only if
$$
\eta (a^{1/p^n}) = a^{1/p^n}
$$
for all $\eta \in N$, if and only if
$$
\QQ(\zpn) (a^{1/p^n}) \subset K_T \;,
$$
if and only if $\QQ(\zpn)(a^{1/p^n})$ is unramified outside $T$, if and only if
$$
v(a) \equiv 0 \mod p^n \hspace{1cm} \forall v \nmid T \;.
$$

\segment{110503b}{}
It remains to verify that the map $\ka$ is equivariant with respect to the action of $\Zpnstar$. To this end, fix $\al \in \Zpnstar$, $a \in E$, $\eta \in \HTn$, and a $\ga \in G_T$ mapping to $\eta$. Then we have
\begin{align*}
\ka_a^\al(\eta) 	&= (\ka_a(\ga \eta \ga \inv))^{\al^{-2}} \\
				&= \left( \frac {\ga \eta \ga \inv (a^{1/p^n}) } { a^{1/p^n} } \right)^{\al^{-2}} \\
				&= \left(  \ga \frac { \eta (\ga\inv a)^{1/p^n} }  { (\ga\inv a)^{1/p^n} }  \right)^{\al^{-2}} \\
				&= \left( \frac { \eta (\al\inv a)^{1/p^n} }  { (\al\inv a)^{1/p^n} } \right)^{\al\inv} \\
				&= \frac { \eta (\al\inv a^{\al\inv})^{1/p^n} }  { (\al\inv a^{\al\inv})^{1/p^n } } \\
				&= \ka_{\al\inv a^{\al\inv}}(\eta) \;,
\end{align*}
indeed.

\segment{110503a}{Proposition}
Let $\Ttilde$ denote the set of primes of $\ZZ[\zpn]$ above $T$. We identify $\Ttilde$ with the set of corresponding valuations of $\QQ(\zpn)$.
Given $\al \in \Zpnstar$ and $b \in (\Zpn)^\Ttilde$, we let $\al$ act on $b$ by
$$
(\al \star b)_v = \al\inv b_{\al\inv(v)} \;.
$$
Then the formula
$$
a \mapsto (v(a))_{v \in \Ttilde}
$$
defines a $\Zpnstar$-equivariant isomorphism
$$
\frac {E} { \QQ(\zpn)^{*p^n} \ZZ[\zpn]^*} \to {(\Zpn)}^\Ttilde \;.
$$

\begin{proof}
By \cite[Corollary 10.5]{Washington}, the $p$-part of the Picard group of $\ZZ[\zpn]$ vanishes, so multiplication by $p^n$ on $\Pic \ZZ[\zpn]$ is an automorphism. Evaluating the short exact sequence of sheaves
$$
0 \to \Oo^* \to \Kk^* \to \Div \to 0 \;,
$$
together with the endomorphism given by multiplication by $p^n$, on $\ZZ[\zpn]$ (and recalling that on an integral scheme, $\Kk^*$ is flasque), we obtain the following diagram
\[\xymatrix	@ C=20pt	{
			&														&									& 0 \ar[d]							&			\\
0 \ar[r]		& \QQ(\zpn)^*/\ZZ[\zpn]^* \ar[r] \ar[d]	_{p^n}				& \Div \ZZ[\zpn] \ar[r] \ar[d]_{p^n}	& \Pic \ZZ[\zpn]	\ar[r] \ar[d]_{p^n}	& 0			\\
0 \ar[r]		& \QQ(\zpn)^*/\ZZ[\zpn]^* \ar[r] \ar[d]						& \Div \ZZ[\zpn] \ar[d] \ar[r]			& \Pic \ZZ[\zpn]	\ar[r] \ar[d]			& 0			\\	
			& \QQ(\zpn)^*/\ZZ[\zpn]^*\QQ(\zpn)^{*p^n} \ar[d]			& (\Zpn)^{|\ZZ[\zpn]|	_0}	\ar[d] 		& 0									&			\\
			& 0														& 0									&									&			}
\]
in which all rows and columns are exact. Here $|\ZZ[\zpn]|_0$ denotes the set of (nonzero) primes of $\ZZ[\zpn]$ (in this notation, $\Div \ZZ[\zpn] = \ZZ^{|\ZZ[\zpn]|_0}$). The snake lemma produces an isomorphism
$$
\QQ(\zpn)^*/\ZZ[\zpn]^*\QQ(\zpn)^{*p^n} \xrightarrow{\cong} (\Zpn)^{|\ZZ[\zpn]|	_0} \;.
$$
It is clear now that the preimage of $(\Zpn)^\Ttilde$ is as stated in the theorem. 

Regarding equivariance, we note that if $K/k$ is a Galois extension, $a\in K$, $\al$ is an automorphism of $K/k$, and $v$ is a place of $K$, then $\al$ induces an isomorphism 
$$
K_v \xrightarrow{\cong} K_{\al(v)} \;,
$$
so $v(\al (a)) = \al(v)(a)$. This completes the proof of the proposition.
\end{proof}

\segment{110505a}{}
The sequence
$$
0 \to \mupn \to \ZZ[\zpn]^* \xrightarrow{p^n} \ZZ[\zpn]^* \to \frac {E} {\QQ(\zpn)^{*p^n}}  \to \frac { E } { \ZZ[\zpn]^*\QQ(\zpn)^{*p^n} } \to 0
$$
is exact.

\medskip
\noindent
This is clear.

\segment{110522a}{}
Summarizing, we have the following diagram of $\Zpnstar$-modules, in which the vertical sequence is exact.
\[\xymatrix{
0 												\ar[d] 					&
																		\\
\frac { \ZZ[\zpn]^* } { \ZZ[\zpn]^{*p^n} } 		\ar[d] 					&
																		\\
\frac {E} {\QQ(\zpn)^{*p^n}} 					\ar[d] \ar[r]^-\cong 		&
\Hom(\HTn, \mupn) 													\\
(\Zpn)^\Ttilde 									\ar[d] 					&
																		\\
0																		}
\]

\medskip
\noindent
We end our study of the terms $\tensor*[^n_n]{E}{^{1,1}_2}$ with a discussion of the structure of 
\[
\frac { \ZZ[\zpn]^* }  { \ZZ[\zpn]^{*p^n} }
\]
as $\Zpnstar$-module.

\segment{110510a}{Proposition}
Denote $\zpn$ by $\ze$ for short. For $a\in \Zpnstar$, let
$$
\xi_a = \ze^\frac{1-a}{2} \frac{1-\ze^a}{1-\ze} \;.
$$
Then we have
\begin{align}
\xi_1 	& = 1 \;, \label{110510d} \\
\intertext{and for each $a \in \Zpnstar$,} 
\xi_a 	& \equiv \xi _{-a} \mod \ZZ[\zpn]^{*p^n} \label{110510e} \;. 
\end{align}
The elements $\xi_a$ of $\ZZ[\zpn]^*/\ZZ[\zpn]^{*p^n}$ parametrized by
$$
a \in \Zpnstar/\langle-1\rangle
$$
are free except for the single relation (\ref{110510d}). If $B$ denotes the $\Zpn$-submodule generated by these elements, then
$$
\ZZ[\zpn]^*/\ZZ[\zpn]^{*p^n} = \mupn \oplus B \;.
$$

\medskip
\noindent
The proof is in paragraphs \ref{110510b}--\ref{110510c}.

\segment{110510b}{}
Equation (\ref{110510d}) is clear. To verify (\ref{110510e}), we carry out the following computation inside $\ZZ[\zpn]^*$:
\begin{align*}
\xi_{-a} 	&= \ze^\frac{1+a}{2} \cdot \frac{\ze^a}{\ze^a} \cdot \frac{1-\ze^{-a}}{1-\ze^a} \cdot \frac{1-\ze^a}{1-\ze} \\
			&= \ze^\frac{1+a}{2} \cdot \frac{-1}{\ze^a} \cdot \frac{1-\ze^a}{1-\ze} \\
			&= -\xi_a \;,
\end{align*}
and note that
$$
-1 = (-1)^{p^n} \equiv 1 \mod \ZZ[\zpn]^{*p^n} \;.
$$

\segment{110510c}{}
Let $U:= \ZZ[\zpn]^*$, let $C^+$ denote the subgroup generated by the elements $\xi_a$, $a \in \Zpnstar$, and let $U^+$ denote the subgroup of $U$ of totally real units. Then by \cite[Theorem 8.2]{Washington}, $C^+$ is a subgroup of $U^+$ of index $h^+$, the class number of the maximal totally real subfield. By \cite[Theorem 4.12]{Washington}, $\mupn \oplus U^+$ has index $1$ or $2$ in $U$. Since $h^+|h$ and $h$ is coprime to $p$, it follows that $\mupn \oplus C^+ \le U$ is a subgroup of finite index coprime to $p$. According to the Dirichlet unit theorem,
$$
U \cong \mupn \oplus \ZZ^{r+s-1}
$$
where $r$ is the number of real places, and $s$ is the number of complex conjugate pairs of complex places. Thus, in our case, $r= 0$ and
\[
s	= \big| \Zpnstar / \langle-1\rangle \big| \;.
\]
It follows that the $\xi_a$ generate a free abelian group of rank $s-1$, and that their image modulo $U^{p^n}$, together with $\mupn$, generates all of $U/U^{p^n}$.  This completes the proof of the proposition.

\segment{110510f}{}
Let $\Zpnstar$ act on $\ZZ[\zpn]^*/\ZZ[\zpn]^{*p^n}$ by $\be \circ a = \be (a)^\be \;.$ This is the action induced by the action defined in paragraph \ref{110502h}, except for having taken the liberty to precompose with the automorphism of $\Zpnstar$ given by $\al \mapsto \al\inv$. We recall that here multiplication by $\be$ on the left refers to the Galois action, while the exponent refers to multiplication inside $\ZZ[\zpn]^*$.

\segment{110510g}{Proposition}
We have
\begin{align}
\be \circ \ze 	& = \ze^{\be^2} 						\label{110510h} \\
\intertext{for any $\ze \in \mupn$, and}
\be \circ \xi_a 	& = \xi_{\be a}^\be \xi_\be^{-\be} \;.	\label{110510i}
\end{align}

\begin{proof}
Equation (\ref{110510h}) is clear. To verify (\ref{110510i}), we compute, focusing on the Galois action:
\begin{align*}
\be (\xi_a) 		&= \ze^{\be \frac{1-a}{2}} \frac {1-\ze^{\be a}} {1-\ze^\be} \\
	& = \ze^{ \frac {1-\be a} {2} - \frac {1-\be} {2}} \cdot \frac {1-\ze^{\be a}} {1-\ze} \cdot \left( \frac {1-\ze^\be} {1-\ze}  \right)\inv 	\\
	& = \xi_{\be a}\xi_\be\inv \;.
	\qedhere
	\end{align*}
\end{proof}

\bibliography{references}

\begin{thebibliography}{DCW}

\bibitem[DCW]{CKTwo}
I.~Dan-Cohen and S.~Wewers.
\newblock Explicit chabauty-kim theory for the thrice punctured line in depth
  two.
\newblock Preprint.

\bibitem[Ser]{corloc}
Jean-Pierre Serre.
\newblock {\em Corps locaux}.
\newblock Hermann, Paris, 1968.
\newblock Deuxi{\`e}me {\'e}dition, Publications de l'Universit{\'e} de
  Nancago, No. VIII.

\bibitem[Sha]{Sharifi}
Romyar Sharifi.
\newblock On a result of {S}oul\'e.
\newblock Unpublished, available at
  http://math.arizona.edu/[tilde]sharifi/soule2.pdf, 2000.

\bibitem[Was]{Washington}
Lawrence~C. Washington.
\newblock {\em Introduction to cyclotomic fields}, volume~83 of {\em Graduate
  Texts in Mathematics}.
\newblock Springer-Verlag, New York, 1982.

\end{thebibliography}

\bibliographystyle{alphanum}

\bigskip

\Small\textsc{I.D.: Fakult\"at f\"ur Mathematik, 
Universit\"at Duisburg-Essen,
Universit\"atsstrasse 2, 
45117 Essen, 
Germany}

\Small\textsc{S.W.: Universit\"at Ulm,
Institut f\"ur Reine Mathematik,
Helmholtzstrasse 18,
89069 Ulm, Germany}

\smallskip

\textit{E-mail address (I.D.):} \texttt{ishaidc@gmail.com}

\smallskip

\textit{E-mail address (S.W.):} \texttt{stefan.wewers@uni-ulm.de}

\end{document}